\documentclass[11pt]{article}
\usepackage{amsmath,amsthm,amsfonts,latexsym,amssymb,enumerate,color}
\usepackage{graphicx}

\def\CC{{\mathbb C}}

\def\GL{{\mathcal L}}
\def\GX{{\mathcal X}}

\def\phi{\varphi}

\def\rank{\mathop{\rm rank}}
\def\re{\mathop{\rm Re}}

\def\RR{{\mathbb R}}
\def\supp{\mathop{\rm supp}\nolimits}

\newtheorem{thm}{Theorem}[section]
\newtheorem{prop}[thm]{Proposition}

\newtheorem{lem}[thm]{Lemma}
\newtheorem{cor}[thm]{Corollary}
\newtheorem{rem}[thm]{Remark}

\newtheorem{ex}[thm]{Example}

\def\beginpf{\begin{proof}}
\def\endpf{\end{proof}}

\def\beq{\begin{equation}}
\def\eeq{\end{equation}}

\begin{document}

\title{Diffusive systems and  weighted Hankel operators}
\author{Aolo Bashar Abusaksaka\thanks{\tt aolo\_besher@yahoo.com}\ \ and Jonathan R. Partington\thanks{\tt j.r.partington@leeds.ac.uk} \\
School of Mathematics, University of Leeds,\\
Leeds LS2 9JT, U.K.}
\date{}

\maketitle

\begin{abstract}
We consider diffusive systems, regarded as  input/output systems with a kernel given as
the Fourier--Borel transform of a measure in the left half-plane. Associated with these
are a family of weighted Hankel integral operators, and we provide conditions for them
to be bounded, Hilbert--Schmidt or nuclear, thereby generalizing results of Widom, Howland
and others.
\end{abstract}



\medskip

\section{Introduction}
In this paper we explore various operator-theoretic
properties associated with linear time-invariant systems, 
beginning with the comparatively simple property of 
BIBO stability and then considering properties of weighted integral operators, including Hankel operators used in $H^\infty$ approximation
(see, e.g. \cite{GCP})
and the Glover operator used in $L^2$ approximation \cite{GLP3}.

The systems we consider will have impulse responses expressible  as Laplace transforms of measures, and
thus may be discussed using the language of diffusive systems in the sense of Monts\'eny.
In \cite{M},
diffusive systems are defined as SISO linear  time-invariant convolution
systems of the form
\[
y(t)= \int_0^t h(t-\tau) u(\tau) \, d \tau,
\]
where the impulse reponse $h$ is the Laplace transform of a signed measure (or more generally a distribution) $\mu$ defined on $(0,\infty)$;
i.e.,
\[
h(t)= \int_0^\infty e^{-\xi t} \, d\mu(\xi) \qquad (t \ge 0).
\]
The associated transfer function is the Stieltjes transform of $\mu$, given by the formula
\[
G(s)=(\GL h)(s)= \int_0^\infty e^{-st} h(t) \, dt = \int_0^\infty \frac{d\mu(\xi)}{s+\xi}
\]
for $s \in \CC_+$, the open right half-plane.
As explained in \cite{M} a diffusive system with measure $\mu$ can be realized in terms of the heat equation
\[
\Psi_t(x,t)=\Psi_{xx}(x,t)+\delta(x) u(t)
\]
with $\Psi(x,0)=0$ $(x \in \RR)$, and
\[
y(t)=\int_{-\infty}^\infty 4\pi^2 x   \Psi(x,t) \, d\mu(4\pi^2 x^2).
\]
Some advantages of diffusive representations are that we may
represent causal convolutions as classical input/output dynamical systems.
This allows the use of a range of PDE techniques.
Diffusive systems are also 
appropriate for modelling long-memory systems, fractional integrators, etc.\\

More recently, in the book \cite{Mbook} and the tutorial article  \cite{CM}, the 
notion of a diffusive system has been generalized. The starting point is now a mapping
$\gamma \in W^{1,\infty}(J;\CC)$, 
the classical Sobolev space of absolutely continuous functions with $\gamma$, $\gamma'$ bounded; here
$J$ is a subset of $\RR$, defining a closed (possibly at $\infty$) contour 
lying in 
a sector in the left-hand complex half-plane $\CC_-$; in this case
we have the expression
\[
h(t)=\frac{1}{2\pi i} \int_\gamma e^{tp} G(p) \, dp =  \int_J e^{\gamma(\xi)t } \mu(\xi) \, d\xi,
\]
where $G=\GL h$ is the transfer function, and
$\mu(\xi)=\dfrac{\gamma'(\xi)}{2\pi i} G(\gamma(\xi))$.\\

In this note we shall work with a more convenient  definition, 
which is also slightly more general. We take
an arbitrary $\sigma$-finite Borel measure $\mu$ on $\CC_-$ satisfying the
condition
\beq\label{eq:mu}
\int_{\CC_-} e^{ts} \, d|\mu|(s) < \infty \qquad \hbox{for all} \quad t>0.
\eeq
This enables us to define $h$ directly as the Fourier--Borel transform of $\mu$, namely,
\beq\label{eq:fb}
h(t) = \int_{\CC_-} e^{tp} \, d\mu(p),
\eeq
in which case we also have the Stieltjes transform formula
\[
G(s) = \int_{\CC_-} \frac{d\mu(p)}{s-p} \qquad \hbox{for} \quad s \in \CC_+.
\]
Since the functionals $f \mapsto f^{(k)}(a)$ can be expressed using Cauchy integrals for any $a \in \CC_-$ and $k =0,1,2,\ldots$ we see that the
impulse responses $t^k e^{-at}$, and hence all finite-dimensional stable systems, can be represented in this way
using measures $\mu$ (rather than requiring distributions).

\begin{rem}
Yet more general definitions in terms of holomorphic distributions can be found in the thesis \cite{aolo}.
For if we let $\GX$ denote the Fr\'echet space of analytic functions $f: \CC_+ \to \CC$
satisfying the condition that
each of the seminorms
\[
\|f\|_n= \max_{0 \le j \le n} \max_{0 \le k \le j+1} \sup_{z \in \CC_-} |(\re z)^k f^{(j)})(z)| 
\]
is finite, then we may define the Fourier--Borel and Stieltjes transforms of distributions in the dual space of $\GX$, since
$\GX$ contains the exponentials $p \mapsto e^{pt}$ for $t>0$ as well as the kernels $p \mapsto 1/(s-p)$ for $s \in \CC_+$.
\end{rem}

\section{Stability and weighted Hankel operators}

 \subsection{Stability}

The following result may be seen as a natural generalization of the result of Montseny \cite[Thm. 4.4]{M}, which applies to measures on $\RR_+$.

\begin{prop}\label{prop:bibo}
Let $h$ be an impulse response given by the diffusive representation (\ref{eq:fb}), where the associated measure $\mu$ satisfies (\ref{eq:mu}).
If in addition $\mu$ satisfies the condition
\beq\label{eq:stable}
\int_{\CC_-} \frac{d|\mu|(p) }{|\re p|} < \infty,
\eeq
then the impulse response $h$ lies in $L^1(0,\infty)$, thus defining a BIBO-stable system.
For positive measures supported on $(-\infty,0)$ condition (\ref{eq:stable}) is necessary and sufficient for BIBO stability.
\end{prop}

\beginpf
We have
\begin{eqnarray*}
\int_0^\infty |h(t)| \, dt &\le& \int_{t=0}^\infty \int_{p \in \CC_-} |e^{tp}| \, d |\mu|(p) \, dt \\
& = &
\int_{p \in \CC_-} \frac{d|\mu|(p) }{|\re p|} < \infty,
\end{eqnarray*}
by Fubini's theorem, and this
implies the BIBO stability. 

In the case that $\mu \ge 0$ and $\supp \mu \subset (-\infty,0)$, we
have equality in the above, i.e.,
\[
\int_0^\infty |h(t)| \, dt = \int_0^\infty h(t) \, dt = \int_{\RR_-} \frac{d\mu(p)}{|p|}.
\]
Hence if 
(\ref{eq:stable}) fails to hold, the system is not BIBO stable (consider the constant input $u(t)=1$).
\endpf

\subsection{Weighted Hankel operators}

Achievable bounds in model reduction are linked to properties of the Hankel operator $\Gamma$, which
we can define on $L^2(0,\infty)$ by
\[
(\Gamma u) (t) = \int_0^\infty h(t+\tau) u(\tau) \, d\tau   .
\]
For finite-dimensional systems it is a finite-rank operator, and its rank is the McMillan degree of the system.
If $h \in L^1$ the operator $\Gamma$ is compact. So, defining its singular values as
\[
\sigma_k(\Gamma) = \inf \{\|\Gamma-T\|: \rank(T)<k \},
\]
we have $\sigma_k \to 0$. For effective $H^\infty$ model reduction by balanced truncation or optimal Hankel-norm
reduction we require $\Gamma$ to be nuclear (see \cite{GCP,GO}); that is, we require $\sum_{k=1}^\infty \sigma_k < \infty$.
Indeed, the optimal $H^\infty$ error $E_k$ for a degree-$k$ approximation is bounded by
\[
\sigma_{k+1} \le E_k \le \sigma_{k+1}+\sigma_{k+2}+\ldots.
\]

For $L^2$ model reduction, the weighted Hankel operator $\Theta$ 
introduced by Glover \cite{GLP3}, and
defined by
\[
(\Theta u)(t)= \frac{1}{\sqrt \pi} \int_0^\infty t^{-1/4} h(t+\tau) \tau^{-1/4} u(\tau) \, d\tau
\]
plays a significant role.
It satisfies $\|\Theta\|_{HS}=\|h\|_{L^2}$, where HS denotes the Hilbert--Schmidt   norm given by
\[
\|\Theta\|_{HS}^2 =  \sum_{k=1}^\infty \sigma_k^2.
\]
Moreover $\rank(\Theta)$ is the McMillan degree of the system (as for Hankel operators), meaning that $L^2$ errors
for degree-$k$ approximation are bounded below by
\[
\left(\sigma_{k+1}^2 + \sigma_{k+2}^2 + \ldots \right)^{1/2}.
\]

In order to study these and similar operators in the same framework, we define for 
measurable $w: (0,\infty) \to (0,\infty)$
the weighted Hankel operator $\Gamma_{h,w}$ on $L^2(0,\infty)$, by
\beq\label{eq:wtho}
(\Gamma_{h,w} u)(t)= \int_0^\infty w(t) h(t+\tau) w(\tau) u(\tau) \, d\tau,
\eeq
which, if bounded, is self-adjoint whenever $h$ is real-valued.

\begin{thm}\label{thm:nuchank}
Let $w$ satisfy the condition 
$\psi_p \in L^2(0,\infty)$ for each $p \in \CC_-$, where
\[
\psi_p(t)=w(t) e^{pt}.
\]
If 
\beq\label{eq:condnuc}
\int_{\CC_-} \|\psi_p\|^2_2 d|\mu|(p) < \infty,
\eeq
then the weighted Hankel operator $\Gamma_{h,w}$ given by \eqref{eq:fb} and \eqref{eq:wtho} is nuclear.
In the case that $\mu \ge 0$ and $\mu$ is supported on $\RR_-$, Condition \eqref{eq:condnuc} is
necessary and sufficient for nuclearity.
\end{thm}

\beginpf
Clearly by using the Hahn--Jordan decomposition of the real and imaginary 
parts of $\mu$  we may suppose without loss of generality that $\mu \ge 0$.
We now adapt a proof of Howland \cite{howland} and define
an operator
$T_0$ by
\[
T_0 u = \int_{\CC_-} \langle u,\psi_p \rangle \psi_p \, d\mu(p) \qquad (u \in L^2(0,\infty)).
\]
We then have that $T_0=T$ and the nuclear norm of $T$ is 
bounded by
\beq\label{eq:bddnuc}
\|T\| \le \int_{\CC_-} \|\psi_p\|_2^2 \, d\mu(p).
\eeq
Finally, if $\mu \ge 0$ and $\mu$ is supported on $\RR_-$, the
elementary operators $u \mapsto \langle u,\psi_p \rangle \psi_p$ are all positive, and so 
equality holds in
\eqref{eq:bddnuc}.
\endpf

The following corollary contains Howland's result on nuclearity of Hankel operators (the case $\alpha=0$), as well as 
a result on the nuclearity of Glover's operators (the case $\alpha=-1/4$).

\begin{cor}\label{cor:talpha}
Suppose that $\alpha > -\frac12$ and let  $w(t)=t^\alpha$ for $t>0$. Then the
weighted Hankel operator $\Gamma_{h,w}$ is nuclear provided that
\beq\label{eq:alnuc}
\int_{\CC_-} \frac{d|\mu|(p)}{|\re p|^{2\alpha+1}} < \infty.
\eeq
In the case that $\mu \ge 0$ is supported on $\RR_-$, condition \eqref{eq:alnuc} is
necessary and sufficient for nuclearity.
\end{cor}
\beginpf
This follows directly from Theorem \ref{thm:nuchank}, noting that
\[
\|\psi_p\|^2_2 = \int_0^\infty t^{2\alpha}e^{2(\re p)t} \, dt = \int_0^\infty \left(\frac{u}{2x}\right)^{2\alpha} e^{-u} \frac{du}{2x},
\]
where $x=-\re p$ and $u=2xt$.
\endpf

Note that nuclearity of the unweighted Hankel operator implies BIBO stability of the associated
linear system \cite{GCP}, so that Corollary \ref{cor:talpha} directly implies
Proposition \ref{prop:bibo}.

\begin{rem}
The example $d\mu(p)=(\sin p)\, dp$ for $p <0$ leads to $h(t)=1/(t^2+1)$ and a nuclear Hankel
operator (as seen from \cite[Thm 2.1]{howland}), showing that for signed measures condition \eqref{eq:alnuc} is not always necessary for nuclearity.
There are further details and examples in \cite[Chap. 3]{aolo}.
\end{rem}

The Hilbert--Schmidt condition is rather easier to test, but we include the following specimen result for completeness.

\begin{prop}
Suppose that $w(t)=t^\alpha$ with $\alpha>-1/2$. Then
$\Gamma_{h,w}$ is Hilbert-Schmidt if and only if
\[
\int_0^\infty u^{4\alpha+1} |h(u)|^2 \, du < \infty.
\]
If  $\mu \ge 0$ is supported on $\RR_-$, then 
this holds if and only if
\[
\int_{\RR_-}
\int_{\RR_-} \frac{d\mu(x) \, d\mu(y)}{|x+y|^{4\alpha+2}} < \infty.
\]
\end{prop}
\beginpf
It is well known (see e.g. \cite[Chap.~2]{davies}) that an integral operator on a space $L^2(X)$,
given by a measurable kernel $K(s,t)$, is Hilbert--Schmidt if and only
$K \in L^2(X \times X)$.
Since
\begin{eqnarray*}
\int_{t=0}^\infty\int_{\tau=0}^\infty w(t)^2 |h(t+\tau)|^2 w(\tau)^2 \, dt \, d\tau
&=& \int_{u=0}^\infty \int_{\tau=0}^u (u-\tau)^{2\alpha} \tau^{2\alpha} |h(u)|^2 \, du \, d\tau\\
&&  \hskip-55pt=\int_{u=0}^\infty\int_{\lambda=0}^1 u^{4\alpha} |h(u)|^2 (1-\lambda)^{2\alpha}\lambda^{2\alpha} u \, du \, d\lambda,
\end{eqnarray*} 
we have the first expression; then, using the formula \eqref{eq:fb} for $h$, we arrive at
\[
C_1\int_0^\infty  u^{4\alpha+1} \int_{\RR_-}\int_{\RR_-} e^{u(x+y)} \, d\mu(x) \, d\mu(y) \, du
=C_2 \int_{\RR_-}\int_{\RR_-}   \frac{d\mu(x) \, d\mu(y)}{|x+y|^{4\alpha+2}},
\]
where $C_1$ and $C_2$ are constants depending only on $\alpha$.
\endpf

A far more difficult question is the boundedness of $\Gamma_{h,w}$. For unweighted Hankel operators,
much is known: for example, the reproducing kernel thesis holds \cite{bonsall}, meaning that
it is sufficient (and clearly also necessary) that $\sup_{s \in \CC_-} \|\Gamma k_s\|/\|k_s\|< \infty$, where
$k_s(t)=e^{st}$ (these act as reproducing kernels in $H^2(\CC_+)$). However, it is not known
whether this result generalises to weighted Hankel operators.

We shall take an approach based on results in \cite{power} for unweighted Hankel operators, combined with very recent results
from \cite{JPP14}
on Carleson embeddings.

\begin{lem}\label{lem:btoc}
Suppose that $\mu \ge 0$ is 
supported on $\RR_-$,
and that $h$ is given by \eqref{eq:fb}.
Let $w$ be a non-negative weight on $(0,\infty)$.
Define $Z_\mu: L^2(0,\infty) \to L^2(\CC_-,\mu)$ by
\[
Z_\mu f(s) = \int_0^\infty w(t) e^{st} f(t)\, dt.
\]
Then $\Gamma_{h,w}$ is bounded if and only if $Z_\mu$ is bounded, and 
this holds if and only if the reversed Laplace transform $\mathcal R$ given by
\[
(\mathcal R f)(s)=\int_0^\infty e^{st} f(t) \, dt
\]
is a bounded operator from
$L^2(0,\infty; dt/w(t)^2)$ into $L^2(\CC_-, \mu)$.
\end{lem}

\beginpf
If $Z_\mu$ is bounded,
we have \[
\langle Z_\mu f, Z_\mu  g \rangle = \int_{s \in \RR_-} \int_0^\infty w(t) e^{st} f(t)\, dt
\int_0^\infty w(\tau) e^{s\tau}   \overline {g(\tau)}\, d\tau \, d\mu(s)
= \langle \Gamma_{h,w} f,   g \rangle,
\]
so that $\Gamma_{h,w}$ is also bounded. Conversely, putting $f=g$, we see that the boundedness of
$\Gamma_{h,w}$ implies the boundedness of $Z_\mu$.
The mapping $f \mapsto fw$ is an isometry between $L^2(0,\infty)$ and $L^2(0,\infty; dt/w(t)^2)$, and so the
last assertion follows.
\endpf

%

The case $w(t)=1$ is due to Widom, and is equivalent to the condition that $\mu$ is a Carleson measure in the 
space $H^2(\CC_-)$, so that $\mu(-x,0)=O(x)$ as $x \to 0$ and $x \to \infty$ (see \cite{power,widom}).
We are able to extend this result to the class of power weights as follows.  

\begin{thm}\label{thm:boundedcm}
Let $w(t)=t^\alpha$ for $\alpha \in \RR$, and let $\mu \ge 0$ be a measure supported on $\RR_-$.\\
(i) If $-1/2 < \alpha < 0$, then $\Gamma_{h,w}$ is bounded if and only there is a $\gamma>0$ such that
$\mu(-2x,-x) \le \gamma x^{1+2\alpha}$ for all $x>0$.\\
(ii) If $\alpha=0$, then $\Gamma_{h,w}$ is bounded if and only if 
there is a $\gamma>0$ such that
if $\mu(-x,-0) \le \gamma x$ for all $x>0$.\\
(iii) If $\alpha > 0$,  then $\Gamma_{h,w}$ is bounded if and only if 
there is a $\gamma>0$ such that
if $\mu(-x,-0) \le \gamma x^{1+2\alpha}$ for all $x>0$.
\end{thm}

\beginpf
(i) By Lemma  \ref{lem:btoc}, boundness of $\Gamma_{h,w}$ is equivalent to
boundedness of the reversed Laplace transform $\mathcal R:L^2(0,\infty; dt/w(t)^2)$ into $L^2(\CC_-, \mu)$.
By \cite[Thm. 3.11]{JPP14} this is equivalent to the condition
$\mu(-2x,-x) \le \gamma x^{1+2\alpha}$.\\
(ii) This is Widom's result \cite{widom}. See also \cite[Thm.~2.5]{power}. It can also be shown
as in (iii) below.\\
(iii) Again, by Lemma  \ref{lem:btoc}, boundness of $\Gamma_{h,w}$ is equivalent to
boundedness of the reversed Laplace transform $\mathcal R:L^2(0,\infty; dt/w(t)^2)$ into $L^2(\CC_-, \mu)$.
We have $1/w(t)^2=t^{-2\alpha}$, and since for $\beta>-1$ we have
\[
\int_0^\infty e^{-2rt} r^\beta \, dr = \int_0^\infty e^{-u} \left( \frac{u}{2t}\right)^\beta\,  \frac{du}{2t} =
\frac{\Gamma(\beta+1)}{2^{\beta+1}t^{\beta+1}},
\]
it follows from \cite[Prop. 2.3]{JPP13}
that with $2\alpha=\beta+1$ the space $L^2(0,\infty; dt/w(t)^2)$ is isomorphic under the Laplace
transform to a Zen space, which in this case is a weighted Bergman space with weight $|x|^{2\alpha-1}\, dx \,dy$.
Then \cite[Theorem 2.4]{JPP13} implies that $\mathcal R$ is bounded if and only if
\[
\mu(-x,0) \le \gamma x\int_0^x r^{2\alpha-1} \, dr = \gamma' x^{1+2\alpha}
\]
for constants $\gamma,\gamma'>0$.
\endpf

\begin{ex}{\rm
(i) Take $d\mu(x)=dx$, Lebesgue measure. 
Thus $h(t)=1/t$. 
This $\mu$ satisfies the condition (ii) of Theorem \ref{thm:boundedcm}, but not Condition (i) with $\alpha=-1/4$.
Therefore the Hilbert--Hankel
operator defined by
\[
\Gamma u(t)= \int_0^\infty \frac{u(\tau)}{t+\tau} \, d\tau
\]
is bounded on $L^2(0,\infty)$ (as is well-known), whereas the corresponding Glover operator
defined by
\[
\Theta u(t)= \int_0^\infty t^{-1/4}\frac{u(\tau)}{t+\tau}\tau^{-1/4} \, d\tau
\]
is unbounded.\\

(ii) Next take $d\mu(x)=|x|^{1/2} \, dx$ so that 
\[
h(t)=\int_0^\infty e^{-tx}\sqrt{x} \, dx = \frac12 \sqrt{\pi} t^{-3/2}.
\]
Now $\mu$ satisfies Condition (i) with $\alpha=-1/4$, but not Condition (ii).
Therefore we conclude that the Hankel operator
defined by
\[
\Gamma u(t)= \int_0^\infty \frac{u(\tau)}{(t+\tau)^{3/2}} \, d\tau
\]
is unbounded  on $ L^2(0,\infty)$ but the Glover operator defined by
\[
\Theta u(t)= \int_0^\infty t^{-1/4}\frac{u(\tau)}{(t+\tau)^{3/2}}\tau^{-1/4} \, d\tau
\]
is bounded. However $\Theta$ is not Hilbert--Schmidt, since $h \not\in L^2(0,\infty)$,
which incidentally provides an example asked for by K.~Glover in conversation.}

\end{ex}

\begin{rem}
Lemma \ref{lem:btoc} and Theorem \ref{thm:boundedcm} have partial extensions to sectorial measures supported on $\CC_-$.
The key observation is that we now have 
\[
\langle Z_\mu f, \overline{Z_\mu g} \rangle= \langle \Gamma_{h,w}f, \overline g\rangle,
\]
so that boundedness of the Laplace--Carleson embedding is sufficient (although possibly not always necessary)
for the boundedness of the integral operator. We leave the details to the interested reader.
\end{rem}

\end{document}